\documentclass[10pt]{amsart}
\usepackage[]{amsmath, amsthm, amsfonts, amssymb}
\usepackage[all]{xy}

 \newcommand {\C} {{\mathbb C}}
 \newcommand {\R} {{\mathbb R}}
 \newcommand {\Z} {{\mathbb Z}}
 \newcommand {\Q} {{\mathbb Q}}

 \newcommand {\HH} {{\mathcal H}}

 \newcommand {\E} {{\mathcal E}}
 \newcommand {\dt} {{\bullet}}
 
 \newcommand {\G} {{\mathcal G}}

 \newcommand {\Gm} {\mathbb{G}_m}

\newcommand{\Qu}{{\mathcal Q}}
\newcommand {\OO}{{\mathcal O}}
\newcommand {\isom}{{\underline{\rm Isom}}}
\newcommand {\fhom}{{\underline{\rm Hom}}}
\newcommand {\sspec}{{\underline{\rm Spec}}}
\newcommand {\shom}{{\underline{\rm Hom}}}
\newcommand {\Sym}{{\rm Sym}}

\newcommand{\qhs}{\Q\text{-HS}}
\newcommand{\sX}{{\mathfrak X}}

\newcommand {\Ind}{\text{Ind-}}
\newcommand {\M} {{\mathfrak M}}
\newcommand{\sS}{{\mathfrak S}}
\newcommand{\sY}{{\mathfrak Y}}
\newcommand{\sF}{{\mathfrak F}}
\newcommand{\sG}{{\mathfrak G}}
\newcommand{\bun}{{\mathfrak Bun}}
\newcommand {\rmbun} {{\rm Bun}^s}
\newcommand {\sC} {{\mathfrak C}}

\newtheorem{thm}[subsection]{Theorem}
\newtheorem{cor}[subsection]{Corollary}
\newtheorem{lemma}[subsection]{Lemma}
\newtheorem{prop}[subsection]{Proposition}
\newtheorem{defn}[subsection]{Definition}

\newtheorem{ex}[subsection]{Example}
\newtheorem{conj}[subsection]{Conjecture}

\renewcommand{\P}{{\mathcal P}}
\newcommand{\lie}{{\rm Lie}}
\newcommand{\tw}[2][\G]{{ }^{#2}#1}

\newcommand{\ch}[1][\G]{{\sf X}(#1)}
\newcommand{\chd}[1][\G]{{\sf X}(#1)^\vee}

\newcommand{\bund}[2][\G]{{\mathfrak Bun}_{#1}^{#2}}

\begin{document}
 \title[The motive of the moduli stack]
{The Motive of the  moduli stack  of $G$-bundles over the universal curve}

 \author{Donu Arapura}
\address{Department of Mathematics\\
Purdue University\\
West Lafayette, IN 47907\\
U.S.A.}
\thanks{First author partially supported by the NSF}
 \author{Ajneet Dhillon}
\address{Department of Mathematics\\
University of Western Ontario\\
London, Ontario\\
Canada}
\thanks{Second author partially supported by NSERC}

\begin{abstract}
We define relative motives in the sense of Andr\'{e}. After
associating a complex in the derived category of motives
to an algebraic stack we study this complex in the
case of the moduli of $G$-bundles varying over the moduli
of curves.
\end{abstract}

\maketitle

For a classical group $G$, 
let $\rmbun_{G,C}$ denote the moduli space of stable $G$-bundles over
a smooth complex projective curve $C$.  As $C$ varies over the space
of  complex genus $g$ curves $M_g(\C)$, the rational cohomology
$H^i(\rmbun_{G,C},\Q)$ fits together to form a variation of  mixed Hodge 
structure.
One of the aims of this work is to understand this 
variation of mixed Hodge structure. The behaviour for large $i$ seems
rather subtle,  and will be left for the future. In this paper, we
determine this structure completely
for  $i$ less than an explicit constant depending on $G$ and $g$.
The main  point is that in this
range, the cohomology of $\rmbun_{G,C}$ agrees with the cohomology of the moduli stack
$\bun_{G,C}$ of all $G$-bundles on $C$ (section \ref{moduli}). Perhaps
contrary to one's first impression, the stack turns out to be the more
accessible of the two objects for this problem.
We show that for all $i$, the variation of mixed Hodge structure
associated to  $H^i(\bun_{G,C})$  can be built out of the tautological
variation of pure Hodge structure associated to $H^1(C,\Q)$ using
standard linear algebra operations and Tate twists (section \ref{stack}); in particular, it
is pure and the Torelli group acts trivially.

The above statements are deduced from finer results at the motivic
level. Motives come in different flavours, and in this paper we present
yet another, which is a relative theory of pure motives 
over a base $\sS$. The base is allowed to be a quotient of a smooth variety defined over a
subfield of $\C$ by a finite group.
When $\sS$ is a point, the theory reduces to 
Andr\'e's \cite{andre}. In general, the category of motives $M_A(\sS)$
in the present sense, over $\sS$,
 forms a semisimple Tannakian
category. There are realizations  of $M_A(\sS)$  to the category
of $\ell$-adic sheaves over $\sS$ and variations of  Hodge structure
over $\sS\times Spec\, \C$.  To every Artin stack $\sX/\C$, we can associate
a well defined  motive $h^i(\sX)$, which maps to the pure Hodge
structure $Gr_WH^i(\sX)$ under the Hodge realization. Similar
statements hold in the relative case. The main theorem
(theorem \ref{thm:atiyah_motivic})
is a motivic version of the Atiyah-Bott isomorphism, which gives a
precise description of the motive of the stack of $G$-bundles over the
universal curve relative to the moduli stack of all curves.

Our thanks to Pramath Sastry and Clarence Wilkerson for
numerous helpful conversations, both virtual and otherwise.
We would also like to thank BIRS in Banff for their support at the initial
stage of this project.

\section{Stacks}

By a stack ${\mathfrak X}$, 
we will mean an algebraic stack in Artin's sense, which is locally of
finite type over some base scheme. We give a working definition which is sufficient for
our needs, but refer to \cite{laumon} for the full story.
  Given a groupoid $(s,t:R\rightrightarrows U,\ldots)$ in the category of 
algebraic spaces over some base,
such that $s,t$ are surjective and smooth (or \'{e}tale in
the special case of Deligne-Mumford stacks), we can associate a stack called
the quotient stack $\sX= [U/R]$, which we can think of  as an
equivalence class in a sense to be explained. Given a surjective
smooth morphism $V\to U$, we get a new groupoid $(R_V=R\times_{U\times U}
(V\times V)\rightrightarrows  V,\ldots)$ by base change. We define
isomorphism for the quotients as the weakest equivalence relation
such that $[U/R]\cong [V/R_V]$ for all such base changes. 
 In the trivial case, where $s$ and $t$
are the identity, we just get back $U$. So stacks
include algebraic spaces.
Stacks form a category (actually
a $2$-category), where a morphism (or more accurately
a $1$-morphism) is given by a morphism between some pair of defining
groupoids.  Given a stack $\sX$, we get a groupoid valued functor
$Y\mapsto \sX(Y) = Hom(Y,\sX)$ on the category of algebraic spaces. This
functor determines $\sX$ and typically gives the most natural
description of it. In cases of interest to us,
the groupoid comes from an action of a group $G$ on $U$; so $R=G\times U$
with $s,t$ being respectively the projection and action maps.
If the quotient $U/G$ exists as an algebraic space, we get a morphism
$[U/G]\to U/G$ which need not be an isomorphism. The quotient $U/G$ is
usually referred to as the coarse moduli space associated to $[U/G]$, and
it can be characterized by the fact that any morphism of $[U/G]$
to an algebraic space would factor through it uniquely.

We list some key examples along with the  functors they represent.

\begin{ex}
Let  $G$ be a group scheme over a field $\kappa$. The classifying
stack $BG = [Spec\, \kappa/G]$. The universal bundle $EG = [G/G] =*$
maps to $BG$.  The quotient scheme $Spec\, \kappa/G$, is trivial,
while $BG$ is not.
$BG(Y)$ is just the groupoid of principal bundles
on $Y$ and isomorphisms between them.
\end{ex}

\begin{ex}
   Let   $\M_g/Spec\, \Q$ 
be the moduli  stack of smooth projective
curves of genus $g$. This is the quotient (stack) of the Hilbert scheme
of tricanonically embedded curves $H_g$ by the appropriate projective linear
group ${\rm PGL}$. We can also realize this as
the quotient stack of the fine moduli space  $M_{g,n}$ of
curves with level $n\ge 3$ structure by $Sp(2g,\Z/n\Z)$.
The quotient of the universal curve  $C_g\to H_g$ by $PGL$
group yields a  a morphism of stacks $\sC_g\to \M_g$.
$\M_g(Y)$ is the groupoid of genus $g$ curves over $Y$ and $\sC_g(Y)$
is the groupoid of curves with a distinguished section. The coarse
moduli space for $\M_g$ is just the moduli space $M_g$.
\end{ex}

\begin{ex}\label{ex:bunG}
Given an algebraic group $G$ over a field $\kappa$ and a $\kappa$-scheme $S$. 
Let $\bun_{G,X/S}$ be the moduli stack of principal $G$-bundles over a
flat family $X/S$, so that $\bun_{G,X/S}(Y)$ is the groupoid of
$G$-bundles over $X\times_SY$. 
There is morphism of stacks $UG\to X\times \bun_G(X)$ which serves
as the universal $G$-bundle. 
When $X=Spec\, \kappa$, this is the classifying stack $BG$ and $UG=EG$.
There is an action of ${\rm PGL}$ on
$\bun_{G,C_g/H_g}$ and by passage to quotients
we obtain
$$\pi:\bun_{G,\sC_g/\M_g}\to \M_g$$
which is the  object of fundamental interest here. 
\end{ex}

The details of the construction of  $\bun_G=\bun_{G,X/S}$ can be found in
\cite{laumon} for vector bundles ($G=GL$), and we can reduce the
general case to this. Choose a faithful representation $G\hookrightarrow GL$, 
then the fibre of the natural map
\[
\rho:\bun_{G}\rightarrow\bun_{\rm GL}
\]
over a principal ${\rm GL}$-bundle
$E$ is the set of  all reductions of the structure group of $E$ to $G$, or in other words
sections of $E/G\to X$. Since $E/G$ is a quasiprojective variety over $X$,
its sections are representable by a subscheme of a Hilbert scheme. This argument
carried out in families shows that the morphism $\rho$ is
representable, which implies that $\bun_G$ is an algebraic stack.

If $X$ is a smooth projective curve over a scheme $S$ and $\G$
is a smooth group scheme over $X$ then $\bun_\G$ is also an algebraic stack.
This requires more work and the details can be found in \cite{behrend}.
\bigskip

A  groupoid  can be extended to a 
simplicial algebraic space  called its nerve
$$\ldots R\times_{s,U,t} R\times_{s,U,t} R\stackrel{\rightarrow}{\rightrightarrows}
 R\rightrightarrows U$$
 When working over $\C$,
the geometric realization of the corresponding analytic object gives a 
topological space $|[U/R]|$ whose weak homotopy type 
 depends only on the underlying stack. 
Thus we have well defined notions
of singular cohomology and fundamental groups for stacks over $\C$.
This can be refined to show that stacks carry natural mixed Hodge
structures on cohomology (of possibly infinite dimension). See
\cite{dhillon, teleman} for details.
More generally, given a stack $F:\sX\to S$ over a complex base scheme.
we can define direct images $R^iF_*\Q$ by the same procedure.
We record the following lemma which is straightforward.

\begin{lemma}\label{lemma:pi1stack}
If $X$ is a connected variety on which a finite group $G$ acts,
and $\sX = [X/G]$. 
There is an  exact sequence of fundamental groups
$$1\to \pi_1(X)\to  \pi_1(\sX)\to G\to 1$$
and  $\pi_1(\sX)$ is
determined by the action of $G$ on the fundamental groupoid of $X$.
\end{lemma}

When $\sX = BG$, the geometric realization of the
nerve is just the usual bar construction for the classifying space.
And for $\sX = [U/G]$, this is nothing but 
the homotopy quotient $(U\times EG)/G$.
When  $\sX=\M_g$, the {\em rational} cohomology is
the same as for the moduli space $M_g$. However, the fundamental groups
are different. When $X=M_{g,n}$, the sequence in the lemma is 
$$1\to \Gamma_{g,n}\to \Gamma_g\to Sp(2g, \Z/n\Z)\to 1$$
where $\Gamma_g$ is the mapping
class group, and  $\Gamma_{g,n}$ is its $n$th congruence subgroup.

\section{Relative Motives}

Relative Chow and homological motives have been constructed by
Denninger-Murre \cite{deninger} and Corti-Hanamura \cite{corti} respectively.
Our goal is to define a relative version of Andr\'e's category
of motives, whose construction we now recall \cite{andre}. 
Let $\kappa$ be a field of characteristic $0$ embeddable into $\C$ with
 algebraic closure   $\bar \kappa$. We will fix this notation for the
remainder of the paper.
Let $SPVar_\kappa$ be the category of smooth
 projective (possibly reducible) varieties over $\kappa$.
Fix  the Weil cohomology $H^*(X)= H^*(X_{et},\Q_\ell)$ for the moment.
 Andr\'e has constructed a $\Q$-subalgebra $A^\dt_{mot}(X)\subset H^{2*}(X)$ of motivated
cycles on $X$. 
A class $\gamma\in A^\dt_{mot}(X)$ if and only
if there exists an object $Y\in SPVar_\kappa$ %$Y\in \V$ 
and algebraic cycles 
$\alpha,\beta$ on $X\times Y$ such that 
$$\gamma = p_*(\alpha\cup *\beta),$$
where $p:X\times Y\to X$ is the projection, and $*$ is the 
Lefschetz involution with respect to
a product polarization \cite{andre}. Given a second Weil cohomology
${H'}^*$ and comparison isomorphisms $H^*\otimes K\cong {H'}^*\otimes K$,
where $K$ is a common overfield of the coefficient fields,
we can identify motivated cycles with respect to $H^*$ and ${H'}^*$.
Thus the initial choice of Weil cohomology is immaterial. Choosing
${H'}^*$ to be Betti i.e. singular cohomology for $X\times_\sigma Spec\, \C$
(for an embedding $\sigma:\kappa\to \C$),
 we find that motivated cycles map to Hodge cycles.
By varying  $\ell$ and $\sigma$, we see that motivated cycles 
map to absolute Hodge cycles.

The category
$M_A$  of Andr\'e motives, can now be constructed
by following the standard procedure of first constructing a category
whose objects are the same as for $SPVar_\kappa$, but with $Hom(X,Y)$ given
correspondences $A_{mot}^{\dim X}(X,Y)$, and then taking the pseudo-abelian
completion and then inverting the Lefschetz motive. 
Or it can defined in one step \`a la Jannsen \cite[sect. 4]{andre}.
The nice feature is that  $M_A$ is   semisimple $\Q$-linear abelian
category;
 this is only conjecturally true for homological motives.
Let $H^*(X)$ denote $\ell$-adic cohomology with its $Gal(\kappa)$-action or
rational singular cohomology of $X\times_\sigma Spec\, \C$, with its
canonical Hodge structure. Then
each smooth projective variety has a functorial motive $h(X)$, such
that the  functor
$X\mapsto H(X)=\oplus_i H^i(X)$ factors through it. This 
 yields a faithful embedding of  $M_A$ into the category of 
$\ell$-adic representations or Hodge structures.
Moreover, since $A_{mot}$ contains K\"unneth projections, we can decompose
 $h(X) = \oplus_ih^i(X)$ in $M_A$ such that the realization of $h^i(X)$ is $H^i(X)$.
We note also that $A_{mot}^i(X) = Hom(\Q(-i), h(X))= Hom(\Q(-i), h^{2i}(X))$.
Define a motive to have weight $i$ if it isomorphic to summand
of $h^j(X)(k)$ for some $X$ and $j,k$ with $i=j-2k$. 
Any motive $T$ decomposes canonically into a direct sum of 
pure of motives  $w^i(T)$ of weight $i$. Under Hodge realization,
$w^i(T)$ corresponds to the maximal sub Hodge structure of weight $i$.

Let $S$ be a geometrically connected smooth variety over $\kappa$ with an action by 
a finite group $G$. Set $\sS = [S/G]$. (It may be helpful to the
following examples in mind: $G=\{1\}$ so $\sS=S$, or $S=M_{g,n}$ and $\sS=\M_g$.)
Let $SPVar_\sS$ be the category of representable smooth projective morphisms
to $\sS$.  Any object of this category $\sF:\sX\to \sS$ can be pulled back to
a $G$-equivariant smooth projective morphism  $f:X\to S$ and
conversely. We will keep this notation throughout this section.
To every  $\sF:\sX\to \sS$ in $SPVar_\sS$, 
we will define a  motive  $h^0(\sS, R^i\sF_*\Q)\in M_A$ %$=M_A(SPVar_\C)$ 
such that its Betti  realization gives $H^0(S\otimes \C,R^if_*\Q)^G$ for
every embedding $\kappa\subset \C$.
Choose a $G$-equivariant nonsingular compactification  $\bar X$ and
a base point $s\in S(\kappa)$. Then $G$ will act on $h^i(\bar X)$.
Set
\begin{eqnarray*}
  h^0(\sS, R^i\sF_*\Q) &= & im[h^i(\bar X)^G\to   h^i(X_s)]\\
&\cong& 
\frac{h^i(\bar X)^G}{\ker(h^i(\bar X)\to  h^i(X_s))\cap h^i(\bar X)^G}
\end{eqnarray*}
and
$$h^0(\sS,\R \sF_*\Q) =\bigoplus_ih^0(\sS,R^i\sF_*\Q) $$
The arguments of \cite[p. 25]{andre} shows that these are independent
of choices. When $S=\sS$,   \cite[2.6]{arapuraL} would in fact  imply the
existence of well defined  higher Leray motives $h^j(S,R^if_*\Q)$.
This construction can be extended to the general case, but we will skip the
details since $h^0$ is sufficient for the needs of this paper.

We define
$$A^i_{mot}(\sX/\sS) = Hom(\Q(-i), h^0(\sS,R^{2i}\sF_*\Q)) $$
The next proposition will give some useful alternative
descriptions, when $\kappa=\C$.
 We note that $R^{j}\sF_*\Q$ can be pulled back to a local
system on  $(S\times EG)/G$, thus $\pi_1(\sS)$ acts on the fibre $H^{j}(X_s)$.

\begin{prop}\label{prop:AXoverS}
Assume $\kappa=\C$ and let $H^*$ denote rational singular cohomology.
With $\bar X$ and $s$ as above, we have
 $$A^i_{mot}(\sX/\sS) = im(A_{mot}^{i}(\bar X)^G\to H^{2i}(X_s))
= A_{mot}^i(X_s)\cap H^{2i}(X_s)^{\pi_1(\sS)}$$ 
\end{prop}

\begin{proof}
For the first part note that
  \begin{eqnarray*}
    A^i_{mot}(\sX/\sS) &=& Hom(\Q(-i), h^0(S,R^{2i}\sF_*\Q))  \\
         &=& Hom(\Q(-i), im[h^{2i}(\bar X)^G\to h^{2i}(X_s)])\\
         &=& im[Hom(\Q(-i),h^{2i}(\bar X))^G \to Hom(\Q(-i),h^{2i}(X_s)) ]\\
         &=& im[A_{mot}^{i}(\bar X)^G\to A_{mot}^{i}(X_s)]\\
         &=& im[A_{mot}^{i}(\bar X)^G\to H^{2i}(X_s)]
  \end{eqnarray*}
The last equality follows from the functoriallity of motivated cycles
under restriction with respect to $X_s\subset \bar X$.

From lemma \ref{lemma:pi1stack}, we see that
$H^{2i}(X_s)^{\pi_1(\sS)}= (H^{2i}(X_s)^{\pi_1(S)})^G$.
Thus we have 
\begin{eqnarray*}
im[A_{mot}^{i}(\bar X)^G\to H^{2i}(X_s)]
&\subseteq&  A_{mot}^i(X_s)\cap (H^{2i}(X_s)^{\pi_1(S)})^G\\
&=& A_{mot}^i(X_s)\cap H^{2i}(X_s)^{\pi_1(\sS)}
\end{eqnarray*}
So it suffices to  prove the reverse inclusion.
Suppose that $\xi_s\in A_{mot}^i(X_s)\cap H^{2i}(X_s)^{\pi_1(\sS)}$.
Since it is both $\pi_1(S)$ and $G$-invariant, it can transported to a $G$-invariant
cycle $\xi_t$ on any fibre $X_t$.
By \cite[thm 0.5]{andre}, this is again motivated.
So we can express $\xi_t = p_*(\alpha_t\cup *\beta_t)$ for algebraic cycles $\alpha_t,\beta_t$
on $X_t\times Y$ for some $Y$. By taking $t$ general, we can (by standard Hilbert scheme
arguments) assume that these
cycles extend to algebraic
cycles $\alpha,\beta$ on $f^{-1}(U)\times Y$ for some Zariski neighbourhood of $t$. Then
by taking closures, we get algebraic cycles
$\bar \alpha,\bar \beta$ on $\bar X\times Y$ such that
$$\xi_s = p_*(\bar \alpha\cup *\bar \beta)|_{X_s} \in  im\,
A_{mot}^i(\bar X) $$
Furthermore by averaging over the group, we can assume these cycles are $G$-invariant.
\end{proof}

By Artin's comparison theorem, we obtain
\begin{cor}
  $$A^i_{mot}(\sX/\sS) 
= A_{mot}^i(X_s)\cap H^{0}(S_{et},R^{2i}f_*\Q_\ell)^G$$ 
\end{cor}

Define the category  $Cor_{A_{mot}}(\sS)$ of relative motivated correspondences 
having as objects smooth projective morphisms, and 
$$Hom(\sX/\sS,\sY/\sS) = \prod_i A_{mot}^{\dim X_i-\dim S}(\sX_i\times_\sS \sY),$$
where $X_i$ are the connected components of $X= \sX\times_\sS S$. 
Composition is defined by  the usual rule \cite{kleiman}.

The category $M_A(\sS)^{st}$  of strict motives over $\sS$ can be constructed
from the category   $Cor_{A_{mot}}(\sS)$  as  above.  
Each smooth projective map $\sF:\sX\to \sS$ gives rise to
a motive $h(\sX/\sS) = h(\sF)$. As before, we have 
$$Hom_{M_A(\sS)}(h(\sX/\sS), h(\sY/\sS)) = Hom_{Cor_{A_{mot}}(\sS)}(\sX/\sS,\sY/\sS)$$

We define realization
functors  from $M_A(\sS)$ to the category of $G$-equivariant polarizable
variations of Hodge structures on $S\times_\sigma Spec\, \C$ by
$$H^i(h(f)(m)) = R^{i+2m}f_*\Q(m) $$
$$H(h(f)(m)) = \bigoplus_i H^i(h(f)(m))$$
We can define realizations to $G$-equivariant $\ell$-adic locally constant
sheaves in a similar fashion.
Given $\sF:\sX\to \sY,\sG:\sY\to \sS$ in $SPVar_\sS$, 
we say that $\sG$ is motivated by $\sF$ if
it lies in the subcategory of  $M_{A}(\sS)$ generated from $\sF$
by taking sums, summands, and products.  
It follows that if $\sG$ is motivated by $\sF$, then
$H(\sG)$ lies in the tensor category generated by 
$H(\sF)$ and Tate structures.

\begin{thm}\label{thm:relmotives}
  The category $M_A(\sS)^{str}$ is a  semisimple Tannakian
 category, and the realization functors
gives exact faithful embeddings of this into
the Abelian categories of polarizable Hodge structures and locally constant
$\ell$-adic sheaves.
\end{thm}

We reduce this to a series of lemmas.

\begin{lemma}\label{lemma:semisimple}
  Let $R$ be a finite dimensional algebra over  $\Q$,
such that it possesses a trace $\tau:R\to \Q$ and an algebra involution
$a\mapsto a'$ such that the bilinear form $\tau(ab')$ is positive definite.
Then $R$ is semisimple.
\end{lemma}

\begin{proof}
  \cite[3.13]{kleiman}.
\end{proof}

\begin{lemma}\label{lemma:involution}
We use the same notation and assumptions as in
proposition~\ref{prop:AXoverS}, in particular that $\kappa=\C$.
Under the inclusion $A_{mot}^*(X/S)\subset H^*(X_s)$,
$A_{mot}^*(X/S)$ is invariant under the  Hodge involution, as defined
in \cite[pp. 10-11]{andre},
with respect to an ample $G$-invariant line bundle (which exists) on $\bar X$.
\end{lemma}

\begin{proof}
Pick an ample line bundle $\mathcal{L}$ on $\bar X$, replace it with $\otimes_{g\in G}\, g^*\mathcal{L}$,
 and equip $\bar X$ and $X_s$ with the associated K\"ahler metrics.
The Hodge involution $*$ is the same as the Hodge star operator (up to a factor and 
complex conjugation). This  is $G$-equivariant since the metric is
invariant. One can check that the Hodge star operator is compatible with
restriction of K\"ahler manifolds. 
Therefore $H^*(X)^{\pi_1(S)}$, which equals $ im(H^*(\bar X)\to H^*(X_s))$
 by the theorem of the fixed part \cite[4.2]{deligne}, is stable under
 $*$, and hence so is $H^*(X)^{\pi_1(\sS)}$.
The invariance  is also true for $A_{mot}^*(X_s)$ by \cite[2.2]{andre}, and thus
 also for the intersection of these spaces. 
\end{proof}

By comparison, we get the same conclusion in general.

\begin{cor}
  The lemma holds for $\ell$-adic cohomology for any field $\kappa$ of
  characteristic zero.
\end{cor}

  \begin{proof}[Proof of theorem \ref{thm:relmotives}]
  To prove that $M=M_A(\sS)$ is Abelian and semisimple,
it is enough by \cite[lemma 2]{jannsen} to prove that
$End_{M}(T)$ is semisimple for each motive $T$.   We can assume that 
$T=h(\sX/\sS)$ for some smooth projective morphism of relative dimension $d$ . 
Then by lemma \ref{lemma:involution}
$$End_{M}(T) = A_{mot}^d(\sX\times_\sS \sX/\sS)\subset H^*(X_s\times X_s)$$
 is stable under the involution
$a' = *a^t*$, where $a^t$ is the transpose
(c.f. \cite[1.3]{kleiman}). 
With the help of  the Hodge index theorem, we see that this algebra
satisfies the
conditions of lemma \ref{lemma:semisimple} (see [loc. cit, p. 381]). Therefore
it is semisimple.

The functor $H$ is exact, since any additive functor between semisimple
categories is exact. The faithfulness follows from
proposition~\ref{prop:AXoverS}.
The tensor structure is induced by fibre product:
$h(\sX/\sS)\otimes h(\sY/\sS) = h(\sX\times_\sS \sY/\sS)$. 
The verification that this
is Tannakian is essentially the same argument
as in \cite[4.3]{andre} and \cite[cor. 2]{jannsen}.
  \end{proof}

Given  a morphism $\mathfrak{T}\to \sS$ of stacks satisfying the above
assumptions, we have a base change
map $M_A(\sS)^{str}\to M_A(\mathfrak{T})^{str}$ for the categories of motives.

\begin{cor}
The base change map is always exact.
Given a point $s\in S$,
the functor $M_A(\sS)^{str}\to M_A(s)^{str} = M_A$ is  
an exact faithful embedding.
\end{cor}

\begin{proof}
As already noted  an additive functor between semisimple categories
is exact.
The composition  $M_A(\sS)^{str}\to  M_A\to \Q_\ell\text{-vect}$ to the category of
$\ell$-adic vector spaces is exact and faithful. 
This also factors as a composition of
$H$ and the exact faithful fibre functor from $\ell$-adic local
systems to $\Q_\ell\text{-vect}$.
\end{proof}

It will be convenient to define a slightly bigger category of relative
motives $M_A(\sS)$. An object consists of a strict motive in $M_A([U/G])$ for
some nonempty $G$-invariant Zariski open $U\subset S$, such that the
underlying $\ell$-adic local system
extends to $S$, for some fixed $\ell$. It is easily seen to be 
equivalent to requiring the extendibility of topological local system
on $U\times_\sigma Spec\, \C$ for some $\sigma$. In particular, this
notion is independent of the choice of $\ell$ or $\sigma$.
 Given two  motives $R,Q$ defined over open sets
$U$ and $V$ respectively, let 
$$Hom_{M(\sS)}(R,Q) = Hom_{M(\sS)^{str}}(R|_{U\cap V}, Q|_{U\cap V})$$
By a theorem  Griffiths \cite[9.5]{griffiths}, the Hodge realization functor $H$
extends to this category  $M_A(\sS)$. The analogue of 
theorem~\ref{thm:relmotives} is easily checked. As before,
objects in the category $M_A(\sS)$ admit a weight grading
$T= \oplus w^i(T)$ such that $w^i(T)$ maps to the maximal sub variation of
Hodge structure of weight $i$ of the realization of $T$.
We note in passing that the construction also works when $S$ is
singular. In this case, we define a motive on $\sS$ as a  strict motive
on a smooth open $G$-invariant set
$U\subset S$ for which the local system extends.
The notion of a variation of Hodge structure can be extended in a
similar fashion.

\section{The Motive of a Stack}

It will be technically convenient to formally adjoin arbitrary
direct sums to $M_A(\kappa)=M_A(Spec\, \kappa)$.
 We  do this by working in the category
$\Ind M_A(\kappa)$  of ind-objects of $M_A(\kappa)$.
This category is again Abelian and semisimple. Moreover
it has arbitrary direct sums:
$$\bigoplus_{i\in I} A_i = 
\varinjlim_{\text{ finite }J} \left(\bigoplus_{j\in J} A_j\right)
$$
In this section, we show how to associate an object $h(\sX)$ in the
derived category $D^+(\Ind M_A(\kappa))$
for any Artin stack. This would lie in $D^+(M_A(\kappa))$ when $\sX$ has finite type.
Note that by semisimplicity, this
decomposes as $h(\sX) \cong \oplus_i \HH^i(h(\sX))[-i]$ where $
\HH^i(h(\sX))$ is the $i$th cohomology. So $h$ can 
viewed as a graded motive.

Consider a pair $(X,D)$, where $X$ is a smooth projective variety
and $D$ is a normal crossings divisor on $X$ defined over $\kappa$. 
Let $D_1,\ldots D_n$ be the irreducible components of $D$. For a subset
$I\subseteq\{1,\ldots,n\}$ define
\[
D_I = \bigcap _{i\in I} D_i.
\]
Note that $D_\emptyset = X$.

 The indexing set $I$ has a natural order inherited from $[n]$ and we set
$I_k= I - \{i_k\}$ for $1\le k \le |I|$. Also define
\[
D^{(l)} = \coprod_{|I|=l} D_I.
\]
The inclusions $D_I \hookrightarrow D_{I_k}$ induce a natural map
\[
\delta_k : D^{(i)} \rightarrow D^{(i-1)}.
\]
Passing to Andre's category of motives $M_A(\kappa)$ we have maps
\[
h(\delta_k): h(D^{(i-1)}) \rightarrow h(D^{(i)}).
\]
Dualizing we obtain
\[
h(\delta_k)^* : h(D^{(i)})(-1) \rightarrow h(D^{(i-1)}).
\]

 We define $h(X,D)$ to be the complex (or its image in the derived category) given by
\[
h(D^{(n)})(-n)\rightarrow \cdots \rightarrow h(D^{(0)})
\]
where the differentials are the alternating sums of the $\delta_k$'s,
and $h(D^{(0)})$ is positioned in degree $0$. The following is immediate.

\begin{prop}
  This construction is contravariantly functorial in the pair
$(X,D)$.
\end{prop}

Fix an embedding of the ground field $\kappa\subset \C$, then 
we can regard $M_A(\kappa)$ as a subcategory of the category of
rational pure Hodge  structures $\qhs$ via
the Hodge  realization functor $H$. To simplify notation, we will
often omit this symbol.
This functor  is exact and and so extends to the derived categories. 
We will use the same symbol for 
the derived functor.

There are  functors
$w^i : \qhs\rightarrow\qhs$
which project onto the weight $i$ piece of the Hodge structure.
These are compatible with the previous $w^i:M_A(\kappa)\to
M_A(\kappa)$ in the sense that they commute with $H$.

\begin{prop}\label{prop:piofhXD}
  Let $(X,D)$ be as above and let $U=X\setminus D$. We denote
by $W$ the weight filtration for the mixed Hodge structure on 
$H^*(U,\Q)$. Then there is a canonical isomorphism 
\[
w^j\HH^i(h(X,D)) \cong Gr_W^j H^{j+i}(U,\Q),
\]
(after realization) where $\HH^i$ denote the $i$th cohomology of the complex.
\end{prop}

\begin{proof}
 \cite[3.2.13]{deligne} yields a spectral sequence
$$E_1^{-m,k+m} = H^{k-m}(D^{(m)})(-m) \Rightarrow H^k(U)$$
with $E_2=E_\infty$. Since $E_\infty$ gives
the weight graded subquotients of the abutment, and
the $E_1$ complex coincides with $w^*\HH^*(h(X,D))$, the result 
is certainly true qualitatively. But we need to calculate the
 precise indices:
\begin{eqnarray*}
w^j\HH^i(h(X,D))  &=& \HH^i(w^j(h(X,D))\\
 &=&  \HH^i(\ldots\to  H^{j-2m}(D^{(m)})(-m)\to \ldots)\\
 &=& E_2^{i,j}\\
 &=& Gr_W^j H^{j+i}(U)
\end{eqnarray*}
In the first line, we use the fact that $w^j$ is exact and thus commutes
with $\HH^i$.
\end{proof}

If we regard $\C_{D^{(i)}}$ as sheaves on $X$, then we form the complex
$$\C_{D^{(n)}}[-2n]\to\ldots \to\C_{D^{(0)}}[0]$$
in the derived category, where the differentials are alternating sums of Gysin maps.
Let $\sigma^k$ denote the stupid filtration \cite[1.4.7]{deligne}.
As a corollary of the proof, we obtain:

\begin{cor}
There is an isomorphism 
 $\HH^i(h(X,D))\otimes \C \cong H^i(\C_{D^{(\dt)}}[\dt])$
under which images under $w^\dt$ on the left correspond to $\sigma$-graded components
on the right.
\end{cor}

It will be important to chose a canonical representative for  this
complex.  Via Poincar\'e residues, we can realize $\C_{D^{(k)}}[-2k]\to\C_{D^{(k-1)}}[2-2k]$  as
the morphism $ \Omega^\dt_{D^{(k)}}[-2k]\to  \Omega^\dt_{D^{(k-1)}}[2-2k]$ in the sense of derived
categories given by
$$
\xymatrix{
    & Cone(W_k[-k]\to \Omega^\dt_{D^{(k)}}[-2k])  &      \\ 
 \Omega^\dt_{D^{(k)}}[-2k]\ar[ru] &   W_{k-1}[1-k]\ar@{=>}[u]\ar[r]&  \Omega^\dt_{D^{(k-1)}}[2-2k] 
}
$$
where $W_\dt = W_\dt\Omega_X^\dt(\log D)$ and the double arrow is a quasi-isomorphism. 
Then by repeated use of mapping  cones (as in \cite[pp 161-162]{gelfand}),
we can use these maps to build a complex of sheaves $K^\dt_{(X,D)}$ on $X$ quasi-isomorphic to $\C_{D^{(\dt)}}[\dt]$,
which is functorial in the pair $(X,D)$.

We can construct a mild generalization of the above. 
Let $(X_\bullet,D_\bullet)$ be a simplicial logarithmic pair defined
over $\kappa$. 
Then complexes $h(X_k,D_k)$ fit into a double complex in $M_A(\kappa)$, where 
the second differentials  are alternating sums of face maps (with signs
suitably adjusted to anticommute). Let $h(X_\dt,D_\dt)$  be the total
complex in $C^+(M_A(\kappa))$. This determines an object in $D^+(M_A(\kappa))$ denoted by
the same symbol. Similarly, we can build a complex of sheaves $K^\dt_{(X,D)}$ on $X_\dt$.

\begin{thm}\label{thm:hXD}
Let $U_\bullet = X_\bullet\setminus D_\bullet$.
  We have 
\[
w^j\HH^i(h(X_\bullet,D_\bullet)) = Gr_W^{i+j} H^j(U_\bullet,\Q).
\]
 (after realization).
\end{thm}

\begin{proof}
  We have a spectral sequence of MHS
  \begin{equation}
    \label{eq:specseq1}
{}_IE_1^{pq} = H^q(U_p)\Rightarrow H^{p+q}(U_\dt)    
  \end{equation}
by \cite[8.3.4]{deligne}. After tensoring by $\C$,
this can be constructed as the spectral sequence for total direct image 
$$Tot(\R\Gamma(\Omega^\dt_{X_0}(\log D_0))\to \R\Gamma(\Omega^\dt_{X_1}(\log D_1))\to\ldots )   $$
associated to the filtration by  skeleta
$$Tot(\ldots 0\to \R\Gamma(\Omega^\dt_{X_p}(\log D_p))\to 
\R\Gamma(\Omega^\dt_{X_{p+1}}(\log D_{p+1}))\ldots )   $$
On the other hand, filtering the double complex 
defining $h(X_\dt, D_\dt)$ by skeleta, yields a spectral sequence
$${}_{II}E_1^{pq} = \HH^q(h(X_p,D_q)) \Rightarrow \HH^{p+q}(h(X_\dt, D_\dt))$$
Therefore after applying $w^j$ and using proposition \ref{prop:piofhXD}, we get
a spectral sequence
$$E_1^{pq}= Gr_W^j H^{q+j}(U_p)\Rightarrow w^j\HH^{p+q}(X_\dt, D_\dt)$$
So it suffices to show that this coincides with the spectral sequence
resulting from applying $Gr^j_W$ to \eqref{eq:specseq1} and shifting. But this
follows from the previous discussion, since we can construct a spectral sequence
using the skeletal filtration on 
$$Tot(K_{(X_\dt, D_\dt)})$$
which maps to both ${}_{I}E_1$ and ${}_{II}E_1$.
\end{proof}

%% Added

\begin{cor}
For any simplicial algebraic space $Y_\dt$ over $\kappa$, with each
$Y_n$ of finite type,
we have an object $h(Y_\dt)\in D^+(M_A(\kappa))$
satisfying
$$w^j\HH^i(h(Y_\bullet)) = Gr_W^{i+j} H^j(Y_\bullet,\Q)$$
 (after realization).
\end{cor}

\begin{proof}
By standard arguments (cf. \cite[8.3.6]{deligne}),
 we can construct a simplicial logarithmic pair $(X_\dt, D_\dt)$
such that $U_\dt=X_\dt-D_\dt$ has the same cohomology as $Y_\dt$.
Moreover, if $(X_\dt',D_\dt')$ is another
such scheme, we can assume without loss of generality that it factors as
$$
\alpha:(X_\bullet', D_\bullet')\rightarrow(X_\bullet,D_\bullet)
$$
Since this induces an isomorphism of mixed Hodge structures on cohomology,
 $h(\alpha)$ must be a quasi-isomorphism by the theorem.
\end{proof}

\begin{cor}
  Let ${\mathfrak X}$ be a stack of finite type over $\kappa$.
Then $h(\sX) = h(Y_\dt)$ gives well define class in
$D^+(M_A(\kappa))$, where $Y_\bullet$ is the nerve of any  
presentation of $\sX$.
\end{cor}

\begin{proof}
This is  well defined, since the mixed Hodge structure on cohomology
depends only on $\sX$ \cite{dhillon, teleman}.
\end{proof}

In view of the above results, it makes sense to define
\begin{equation}\label{eq:hi}
h^i(\sX) = \bigoplus_j w^j \HH^{i-j}(h(\sX)))  
\end{equation}
for a stack or simplicial space. For a smooth projective variety,
this agrees with the previous meaning. In general,
 under Hodge realization,
$h^i(\sX)$ would map to $Gr_WH^i(\sX)=\oplus_jGr_W^j H^i(\sX)$.

We can refine the construction in the following ways:
\begin{enumerate}
\item[(H1)]   For  any stack $\sX$ (locally of finite type as always), we
get a class $h(\sX)\in D^+(\Ind M_A(\kappa))$, such that $h^i(\sX)$
defined as in  (\ref{eq:hi}) maps to the infinite dimensional Hodge
structure $Gr_WH^i(\sX)$.
\item[(H2)]  If $f:X\to S$ is a smooth projective
morphism, and $D\subset X$ a relative normal crossing divisor,
we can define a complex $h((X,D)/S)\in C^+(M_A(S)^{str})$ (and hence a
class in  its derived category) such
that
$$w^j\HH^i(h(X,D)/S) = Gr^j_W R^{j+i}g_*\Q$$
 (after realization) where $g:U\to S$ is the restriction to the complement.
Or equivalently, $h^i((X,D)/S)= Gr_W R^ig_*\Q$ with the above
convention (\ref{eq:hi}). The construction is compatible with base
change.

\item[(H3)] When a finite group $G$ acts equivariantly on $f:X\to S$ and
  $D$ as above, the previous class descends to an element of $ C^+(M_A([S/G])^{str})$.

\item[(H4)] Let $\sS=[S/G]$ be the quotient of a smooth variety by a
  finite group. Call a morphism $\sF:\sX\to \sS$ cohomologically locally
  constant (and finite)
  if all the direct images $R^i(\sF\times_\sS S)_*\Q$ are  locally constant (with
  finite dimensional stalks). Then to
  any  cohomologically locally constant morphism, we can construct
motives $h^i(\sX/\sS)$ in $\Ind M_A(\sS)$ (or  $M_A(\sS)$
  assuming finiteness)  compatible with base change.
\end{enumerate}

Items (H1)-(H3) are straight forward modifications of the previous
construction, so the details will be omitted. However, we will say a
few words about (H4).
Given a cohomologically locally constant  finite morphism $\sF:\sX\to
\sS$, we can find a $G$-equivariant simplicial space
$f_\dt:Y_\dt\to S$, such that  $R^i(\sF\times_\sS S)_*\Q=R^if_{\dt*}\Q$ for
all $i$. We now fix $i$. Then $R^if_{\dt*}\Q$ depends on the $(i+1)$ skeleton which is a
finite diagram. By resolution of singularities applied to the generic
fibre and descent theory \cite{deligne}, we see that there exists   a
nonempty $G$-invariant open set $U\subset
S$, and  a $G$-equivariant  $(i+1)$-truncated simplicial relative
logarithmic pair $g_\dt:(X_\dt, D_\dt) \to U$ such that
 $R^i(\sF\times_\sS S)_*\Q|_U=R^ig_{\dt*}\Q$. Set 
$$h^i(\sX/\sS) =  \bigoplus_j w^j \HH^{i-j}(Tot(h(X_\dt,D_\dt)/U))\in M_A([U/G])^{str}$$
By arguing as in the proof of theorem~\ref{thm:hXD}, we can see that
this maps to $Gr_WR^i(\sF\times_\sS S)_*\Q$ under Hodge realization. Hence it
is  independent of choices and extends to $M_A(\sS)$.

\section{Motive of the moduli stack}
Throughout this section, $G$ denotes a split  semisimple
group over $\kappa$. Fix an embedding $\kappa\subset \C$.
By base change we get a complex group $G$ and stack $BG$ over $\C$,
which  will be denoted by the same symbols when no confusion is likely.
As noted earlier, the space $|BG|$ associated to $BG$
is the classifying space in the usual sense, and we will usually
write $BG$ for both objects to simplify notation.

\subsection{Cohomology of the classifying space}

We recall the description of cohomology of $BG$ and associated spaces.
(The calculations are unchanged if $G$ is replaced by a maximal
compact. After doing so, proofs can be found in \cite{borel, whitehead}.)
We have an  isomorphism  
$$H^*(BG,\Q) \cong H^*(BT,\Q)^W = \Q[x_1,\ldots x_n]^W,$$
where $T=\Gm^n$ is a maximal torus with Weyl group $W$ and $x_i$ is
the first Chern class of universal line bundle on the $i$th factor. 
The right hand side  is a polynomial ring in the elementary $W$-invariant polynomials of the $x_i$.
These are Chern classes of the universal bundle $EG$. Let $2n_i $
denote the degrees of these Chern classes, and let
$N= \bigoplus N_i$ denote the span of these  classes. Since $G$ is
semisimple, these numbers are greater than $2$. These Chern classes define a
map
$$BG\to \prod K(\Z,2n_i)$$
to a product of Eilenberg-Maclane spaces which induces a {\em
  rational} homotopy equivalence. If we identify
$G$ with the based loop space $\Omega BG= Map^*(S^1,BG)$, then we get
a rational homotopy equivalence with $\prod  K(\Z,2n_i-1)$. It follows
that $H^*(G,\Q)$ is an exterior algebra on $N[1]$ (where $N[i]_n = N_{n+i}$).
This can be seen from a different point of view by applying a theorem
of Hopf, then $N[1]$ corresponds to the space
of primitive elements for the Hopf algebra structure. In more explicit
terms, $N[1]$ can be identified with subspace of $H^*(G)$ by taking
the image of $N$ under the composition
\[
 H^*(BG)\stackrel{e^*}{\longrightarrow}
 H^*(G\times S^1)\stackrel{\int_{s^1}}{\longrightarrow} H^{*-1}(G),
\]
of the pullback along evaluation $e:G\times S^1\to BG$ and slant product.
The  loop space  $\Omega G$ which is homotopic to $\Omega^2 BG$ has $|\pi_1(G)|$  connected 
components, each  rationally equivalent
to $\prod  K(\Z,2n_i-2)$. The  cohomology of each component is a symmetric algebra on
$N[2]$. These generators can be obtained from $N$ by the above procedure.

\begin{prop}
  Let $G$ be a split connected semisimple group then 
$h(BG)$ is quasi-isomorphic to a direct sum of translates of Tate motives.
\end{prop}

\begin{proof}
  Let $T$ be a maximal torus inside $G$ and $W$ 
the corresponding Weyl group. The bar construction gives
a model of $BT$ as a simplicial scheme $T_\dt$ for which
each $T_n$ is a union of   a product of $\Gm$'s. It follows
immediately that  $h(BT)$ is a direct sum of translated Tate motives.
Since the cohomology ring of $BT$ is a polynomial ring in $\dim T$
variables,
it follows that there is a quasi-isomorphism
\[
h(BT)\cong \bigoplus_{i=0}^\infty \Q(-i)^{r(i)}[-2i]
\]
where $r(i)=\binom{\dim T+i-1}{i}$.
The Weyl group acts on the right
in a way that is compatible with the action on $h(BT)$.
The natural map $BT\rightarrow BG$ induces an isomorphism
\[
h(BG)\cong (\bigoplus_{i=0}^\infty \Q(-i)^{r(i)}[-2i])^W
\]
as  $H^*(BG) = H^*(BT)^W$.
\end{proof}\subsection{Cohomology of $\bun_G$}

Fix a smooth projective curve $C$ of genus $g$ over $\C$.
Atiyah and Bott  \cite{atiyah} described the cohomology ring of the
mapping space $Map(C, BG)$
and Teleman  \cite{teleman} showed that this space can be identified
with $\bun_G = \bun_{G,C}$. We review these results in a form that is
convenient for us.

In general, $\bun_G$ has $|\pi_1(G)|$ connected components,
Let $\bun^c_G$ denote one of these.
The universal bundle over $C\times\bun^c_G$ produces
a classifying morphism to $BG$. Hence there is a pullback
map
\[
H^*(BG)\rightarrow H^*(C)\otimes H^*(\bun^c_G).
\]
This can be  transposed to obtain
 $$\tau_i: H^*(BG)\otimes H^i(C)^*\to H^{*-i}(\bun^c_G). $$
Set 
$$
\alpha = \tau_0:H^*(BG)\rightarrow H^*(\bun^c_G). 
$$
The maps $\tau_1$ and $\tau_2$ induce algebra homomorphisms
$$
\beta :\bigwedge(N[-1]\otimes H^1(C)^*) \rightarrow H^{*}(\bun^c_G).
$$ 
$$\gamma:Sym(N[-2](1))\rightarrow H^{*}(\bun^c_G).$$

\begin{thm}[Atiyah-Bott, Teleman]\label{thm:atiyah}
  $\alpha\otimes\beta\otimes\gamma$ is an isomorphism of mixed Hodge structures.
\end{thm}

\begin{proof}
This is essentially contained in \cite[section 2]{atiyah} and
\cite[p. 24]{teleman}, but we outline the main points since some
details are only implicit.    Since the map is a morphism of
mixed Hodge structures, it suffices to prove that it is an isomorphism
of vector spaces.  The  Poincar\'e
series of the  domain of  $\alpha\otimes\beta\otimes\gamma$  is  easily computed to obtain
$$
\prod_i\frac{(1+t^{2n_i-1})^{2g}}{(1-t^{2n_i})(1-t^{2n_i-2})}
$$
We can check that we get the same series for each $\bun_G^c$ by
using Thom's theorem (cf \cite[pp 540-541]{atiyah} and
\cite{thom:57}). Note in
particular, that Thom shows that the Poincare polynomial is independent of the
choice of connected component. Thus the Poincar\'e series for $\bun_G$
is $|\pi_1(G)|$ times the above series.

We have a  cofibration
$$\bigvee_{2g} S^1\to C\to S^2$$
which gives a fibration of base point preserving mapping spaces
$$\Omega^2BG=Map^*(S^2, BG) \stackrel{c}{\to} 
Map^*(C,BG)\stackrel{b}{\to} Map^*(\bigvee_{2g} S^1, BG)\sim
\prod_{2g}G$$
where $\sim$ denotes homotopy equivalence.
This together with the  fibration
$$Map^*(C,BG)\to Map(C,BG)\sim \bun_G\stackrel{a}{\to} BG$$
yields a ``3 dimensional spectral sequence''
$$
E_2^{pqr}=  H^p(BG)\otimes H^q( G^{2g})\otimes H^r(\Omega G) \Rightarrow 
H^{p+q+r}(\bun_G)\\
$$
Note that the sum of terms on the left is just a sum of $|\pi_1(G)|$ copies the domain of
$\alpha\otimes\beta\otimes \gamma$. So the equality  of
Poincar\'e series forces $E_2= E_\infty$, and thus we have an isomorphism.
\end{proof}

\begin{cor}
  The mixed Hodge structure on $\bun_G$ is pure, i.e. a direct sum of pure
  Hodge structures.
\end{cor}

\subsection{Motive of the moduli stack}
\label{stack}

Let $C\to S$ be a family of genus $g$ curves. Then we have a $S$-stacks
$BG\times S$ and $\bun_G = \bun_{G,C/S}$. We can apply theorem
\ref{thm:atiyah} fibrewise  to conclude that $\bun_G/S$ is
cohomologically locally constant and finite. Consequently  the motive
$h(\bun_G/S)\in M_A(S)$ is defined.
The universal bundle gives  a morphism
$$h^*(BG\times S/S)\rightarrow h^*(C/S)\otimes h^*(\bun_G/S)$$
as above. It is also clear after passing to the simplicial model,
that we can form the transpose
 $$\tau_i: h^*(BG\times S/S)\otimes h^i(C/S)^*\to h^{*-i}(\bun_G/S) $$
and maps $\alpha,\beta,\gamma$ as we did earlier (where we use the
tensor structure on $M_A(S)$ to define  exterior and symmetric
powers). With this set up, we get as a corollary to theorem
\ref{thm:atiyah}

\begin{cor}
  $\alpha\otimes\beta\otimes\gamma$ is an isomorphism of motives 
\end{cor}

When finite group $\Gamma$ acts on $S$ and the family $C/S$, the above
isomorphism descends to  $M_A([S/\Gamma])$. Applying this to the universal
curve $\pi:\sC_g\to \M_g$ (extended to $Spec\, \kappa$) yields

\begin{thm}\label{thm:atiyah_motivic}
The motive of $\bun_{G,\sC_g}/ \M_g$ is isomorphic to 
  $$h^*(BG\times \M_g/\M_g)\otimes h^*(\bigwedge(N[-1]\otimes
  h^1(\sC_g/M_g)^*))\otimes  h^*(Sym(N[-2](1)) $$
\end{thm}

\begin{cor}
The same isomorphism holds for the variation of Hodge structure
associated to    $\bun_{G,\sC_g}/ \M_g$, and in particular for its
monodromy representation.
\end{cor}

Recall that the Torelli group is the kernel of the monodromy
representation 
$$\Gamma_g\to Sp(2g,\Z)$$
associated to $R^1\pi_*\Z$. As a subcorollary, we see that the
action of Torelli group on cohomology of $\bun_G$ is trivial.

\section{Comparison with the Moduli Space}

Fix a reductive group $G$  and a smooth projective genus $g$
curve $C$, both defined over $\kappa$.

\subsection{(Semi)-Stability for $\G$-bundles}

By a principal $G$-bundle or simply $G$-bundle over $C$, we will mean a scheme $P\to C$,
with a right $G$-action, which is {\'e}tale locally a product.
To every    $G$-bundle $P$ over a curve $C$,
we can form the smooth affine group scheme $\G = Aut(P) =
P\times_{G,Ad} G$. This is reductive since $G$ is.
$\G$ will carry  all the information we need, and it
is technically more convenient to work with it.
We define the \emph{degree} of a smooth affine 
group scheme $\G$ over $C$ to be the degree of the vector
bundle $\lie(\G)$ over $C$. It is denoted by $\deg\G$

The following fact is very useful :

\begin{lemma}
\label{l:trick}
 Let $\G$ be a reductive group scheme over $C$.
  There is a finite \'{e}tale cover $f:Y\rightarrow C$ such that
$f^*\G$ is an inner form.
\end{lemma}

\begin{proof}
  We make use of the notations of \cite{SGA3}. Let $\G_0$ be
the constant reductive group scheme over $C$ having the
same type as $\G$.
Being an inner form
means that the scheme ${\rm Isomext}(\G,\G_0)$ has a section over $C$.
By \cite[XXIV, theorem 1.3]{SGA3} and  by \cite[XXII, corollary 2.3]{SGA3} $G$
is quasi-isotrivial and hence so is ${\rm Isomext}(\G,\G_0)$. This implies
by \cite[X, corollay 5.4]{SGA3} that ${\rm Isomext}(\G,\G_0)$ is \'{e}tale
and finite over $C$. So we take $Y$ to be one of these components
and the section is the tautological section.
\end{proof}

\begin{cor}
If $\G$ is a reductive group scheme over $C$ then
$\deg\G=0$.
\end{cor}

\begin{proof}
  By the above we may assume that $\G$ is an
inner form. If $G_0$ is the constant reductive group
scheme of the same type as $\G$ then the adjoint action
of $G_0$ on its Lie algebra factors through ${\rm SL}(\lie(G_0))$.
\end{proof}

\begin{defn}[Behrend]
  A reductive group scheme $\G/C$  is said to be
\emph{stable} (resp. \emph{semistable}) if for every 
proper parabolic subgroup $\P$ of $\G$ we have
$\deg \P<0$ (resp. $\deg \P \le 0$). The degree of the
largest parabolic subgroup of $\G$ is called the degree
of instability of $\G$ and denoted $\deg_i \G$. 
\end{defn}

Let $E$ be a $\G$-bundle. We denote by $\tw{E}$
the associated inner form
\[
\tw{E}=E\times_{\G,Ad}\G.
\]
We say that $E$ is (semi)-stable if $\tw{E}$ is
and define $\deg_i \G = \deg_i\tw{E}$.

In the case of a constant group scheme $\G=G\times_\kappa C$ this
definition of (semi)-stability is equivalent to the usual
one in \cite{ramanathan}. This equivalence is by the remarks 
in the last paragraph on page 304 of \cite{behrend:95}. Basically
there is a bijection between parabolic subgroups of $\tw{E}$
and reductions of structure group of $E$ to a parabolic subgroup
of $G$. This bijection underlies the equivalence of the
two definitions.

\subsection{The Bounds}

Denote by $\ch=Hom(\G,\Gm)$ the group of characters of $\G$ and by
$\chd={\rm Hom}(\ch,\Z)$ its dual. 

Given a $\G$-bundle $P$
its degree is defined to be the element of $\chd$ 
defined by
\[
\chi\mapsto\deg(P\times_\G{\mathbb G}_m),
\]
where $\chi\in\ch$. Notice that $P\times_\G{\mathbb G}_m$ can be
viewed as a line bundle on our curve.
For $\alpha\in\chd$ we denote by $\bund{\alpha}$
the component of $\bund{}$ consisting parameterizing bundles of degree
$\alpha$. It is a union of connected components of $\bund{}$.

Recall that the stack $\bund{\alpha}$ has dimension $\dim_C \G(g-1)$
where $\dim_C$ is the relative dimension over $C$. For a parabolic subgroup
$\P\subseteq \G$ and a character $\beta\in\chd[\P]$ such that
the universal bundle on $\bund[\P]{\beta}$ has non-negative degree we have 
\[
\dim\bund[\P]{\beta}\le \dim_C\P(g-1).
\]
This follows directly from \cite[proposition 8.1.7]{behrend}.
See also \cite[proposition 5.8]{behrend:06}. Let ${\mathfrak C}$ be the complement
of the stable locus in $\bund{\alpha}$. Observe that, by the definition
of stability, ${\mathfrak C}$ is the union of the images
\[
\bund[\P]{\beta}\rightarrow \bund{\alpha}
\]
where the universal bundle on $\bund[\P]{\beta}$ has non-negative degree
and the image of $\beta$ under
\[
\chd[\P]\rightarrow \chd
\]
is $\alpha$.
It follows from the above that
${\mathfrak C}$ has codimension at least
\begin{eqnarray*}
 & d_\G= {\rm min} (\dim_C\G - \dim_C\P)(g-1) \\
 & = {\rm min}(\dim_C R_u(\P) )(g-1).
\end{eqnarray*}
In the above equation,  ``${\rm min}$'' runs over all proper parabolic subgroups
of $\G$. 

In the case where $\G=C\times_k G$ with $G$ being a split reductive
group, one can interpret the above bound on the codimension in terms
of the root datum of $G$.

First let us recall some facts about parabolic subgroups $P\subseteq G$.
Each such subgroup contains a Borel subgroup $B$ and a maximal torus
$T\subseteq B$. This data determines a set of roots 
$R\subseteq\ch(T)\otimes\Q$, a set of positive roots $R_+$ and a basis
$\Delta$ of $R$.

Let $I\subseteq\Delta$ and let $R_I$ be the set of roots that 
are linear combinations of roots in $I$. Let $W_I$ be the subgroup
of the Weyl group generated by the reflections $S_\alpha$ with $\alpha\in I$.
If 
\[
P_I = \cup_{w\in W_I} B\tilde{w}B\quad
(\tilde{w}\in N_{G}(T)\text{  a representative for } w)
\]
then $P_I$ is a parabolic subgroup of $G$, see 
\cite[Theorem 8.4.3]{springer:98}. Furthermore, there is a 
$J\subseteq\Delta$ so that $P_J= P$. If $\Lambda_I = R_+\setminus R_I$
then 
\[
\dim R_u(P_I) = |\Lambda_I|.
\]
Putting this all together it follows that the bound on the codimension of $C$ 
above is the same as
\[
d_{G} = \min_I (|\Lambda_I|(g-1))
\]
where $I$ runs over all sets of the form $I=\Delta\setminus\{\alpha\}$ for
some root $\alpha\in\Delta$.

We will now study this minimum for the standard families of Chevalley
groups $A_n,B_n,C_n$ and $D_n$.

\begin{prop}\label{prop:dforABCD}
The  minimum values for $d_{G}$ for the standard families of Chevalley
groups $A_n,B_n,C_n$ and $D_n$ are given by the following table\\

\begin{tabular}{|c||c|c|c|c|} \hline
 $G$ &$A_n$ & $B_n\> (n\ge 2)$ & $C_n\> (n\ge2)$ & $D_n\> (n\ge 3)$ \\ \hline
$d_G$ & $n(g-1)$ & $2(n-1)(g-1)$ & $2(n-1)((g-1)$ & $2(n-1)(g-1)$ \\ \hline
\end{tabular}
\end{prop}

\begin{proof}
We do this  case by case.

\underline{$A_n$} This root system can be identified with a subset of the 
hyperplane $H:\sum x_i=0$ in $\R^{n+1}$. If we write $e_i$ for the standard
basis of $\R^{n+1}$ then the roots are $e_i - e_j$, for $i,j$ distinct.
We have the basis
\[
\Delta = \{e_i -e_{i+1} | 1\le i \le n\}.
\]
Hence
\[
R_+ = \{ e_i - e_j|i<j|\}.
\]
A straightforward calculation shows that $d_{G} = n(g-1)$
in this case. Note that this agrees with the bound in
\cite{dhillon} as ${\rm SL}_n$ is of type $A_{n-1}$.

\underline{$B_n$} This root system can be identified with the
subset of $\R^n$ consisting of
\[
\pm e_i\quad\text{and}\quad \pm e_i \pm e_j \quad (i\not= j).
\]
We can take
\[
\Delta = \{\alpha_1,\ldots,\alpha_n\}
\]
where $\alpha_i = e_i - e_{i+1}$ for $i<n$ and $\alpha_n = e_n$.
The positive roots are then :
\begin{eqnarray*}
 \sum_{k=i}^n \alpha_k,\quad &\text{for }1\le i \le n \\
 \sum_{i\le k <j} \alpha_k \quad &\text{for }1\le i<j\le n \\
 \sum_{i\le k<j} \alpha_k + 2(\sum_{j\le k\le n} \alpha_k) \quad &
 \text{for }1\le i<j\le n
\end{eqnarray*}
Let $I=\Delta\setminus\{\alpha_m\}$. The dimension of the unipotent
radical of $P_I$ is then 
\[
m + m(n-m) + \frac{m}{2}(2n-m-1).
\] 
The smallest this can be is $2n-2$, noting that $n\ge 2$.

\underline{$C_n$} This root system can be identified with the subset of
$\R^n$ consisting of
\[
\pm 2e_i\quad \pm e_i\pm e_j\quad i\not= j
\]
We can take
\[
\Delta = \{\alpha_1,\ldots,\alpha_n\}
\]
where $\alpha_i = e_i-e_{i+1}$ for $i<n$ and $\alpha_n = 2e_n$. 
The positive roots are then :
\begin{eqnarray*}
 \sum_{i\le k<j} \alpha_k,\quad &\text{for }1\le i<j \le n \\
 \sum_{i\le k <j} \alpha_k+2\sum_{j\le k<n}\alpha_k+\alpha_n 
      \quad &\text{for }1\le i<j\le n \\
 2\sum_{i\le k<n} \alpha_k + \alpha_n \quad &
 \text{for }1\le i \le n
\end{eqnarray*}
Let $I=\Delta\setminus\{\alpha_m\}$. The dimension of the unipotent
radical of $P_I$ is then 
\[
m(n-m)+ \frac{m}{2}(2n-m-1) + m.
\]
The smallest this can be is $2n-2$, noting that $n\ge 2$.

\underline{$D_n$}This root system can be identified with the subset of
$\R^n$ consisting of
\[
 \pm e_i\pm e_j\quad i\not= j
\]
We can take
\[
\Delta = \{\alpha_1,\ldots,\alpha_n\}
\]
where $\alpha_i = e_i-e_{i+1}$ for $i<n$ and $\alpha_n = e_{n-1}+e_n$. 
The positive roots are then :
\begin{eqnarray*}
 \sum_{i< k<j} \alpha_k,\quad &\text{for }1\le i<j \le n \\
 \sum_{i\le k \le n} \alpha_k
      \quad &\text{for }1\le i < n \\
 \sum_{i\le k<j} \alpha_k + 2\sum_{j\le k\le l-1} \alpha_k +\alpha_{l-1} +
\alpha_l\quad &
 \text{for }1\le i <j < n
\end{eqnarray*}
Let $I=\Delta\setminus\{\alpha_m\}$. The dimension of the unipotent
radical of $P_I$ is then 
\[
(m-1)(n-m)  + m + \frac{m}{2}(2n-m-3).
\]
The smallest this can be is $2n-2$, noting that $n\ge 3$.
\end{proof}

\subsection{The  moduli space}
\label{moduli}

Let $G$ be a split semisimple group over $\kappa$. In this case
the character group $\ch[C\times_\kappa G]$ is trivial. Let $S$
be a scheme over $\kappa$ with an action of an algebraic group $K$.
Let $f:C\to S$ be a smooth projective curve over $S$ for which the action
of $K$ lifts.
Denote by $\bun^{s}_{G}$ the open substack of
$\bun_G=\bun_{G,C/S}$ parameterizing stable bundles. 
Ramanathan  \cite{ramanathan} had originally constructed a coarse moduli space
$\rmbun_G$ for $\bun^{s}_G$, when $S= Spec\, \kappa$.
 We will review its construction in this relative setting, following
 \cite{schmitt:02}.

\begin{thm}

\begin{enumerate}
\item[]
\item[(i)] There is a coarse moduli scheme ${\rmbun}_G$ for the
stack $\bun_G^s$. 

\item[(ii)] The action of $K$ lifts to $\rmbun_G$ and the natural map
$\bun_G^s\rightarrow\rmbun_G$ is equivariant for the $K$-action.
\end{enumerate}
\end{thm}

\begin{proof}

  Fix a faithful representation
 \[
 \rho:G\rightarrow{\rm GL}(V).
 \]
 To give a principal $G$-bundle on a scheme $X$ is the same as giving a
 ${\rm GL}(V)$-bundle $E$ plus a reduction of structure group
 of $E$ to $G$, in other words a section of $E/G$. 
If $\E$
is the vector bundle associated to $E$, then the reduction can
be encoded as a 
homomorphism of sheaves of algebras
\[
\sigma:\Sym^*(\E^\vee\otimes V)^G\rightarrow \OO_X.
\]
such that the induced section of
$$\fhom(V\otimes \OO_X, \E)/G$$
lifts (locally) to a section of $\isom(V\otimes \OO_X, \E)\cong E$.

Note that if a reduction of structure group exists then
$\E$ must have degree $0$. Fix a relatively ample
divisor $D$ on $C$ and set $P(t)=r\deg D t + r(1-g)$
where $r= \dim V$. The collection of vector bundles $\E$
with Hilbert polynomial $P$ that admit reductions to stable
principal $G$-bundles is a bounded family, see \cite[3.2]{schmitt:02}.
So we can find an $N$ so that for any $n\ge N$ and any bundle
in the family
\[
R^if_*(\E(nD)) = 0 \quad\text{and}\quad
\E(nD) \text{ is generated by global sections.}
\]
Let $W$ be a vector space of dimension $P(n)$.
Then  our bounded family is  parameterized by
an open subscheme  $Q$ of the Quot scheme over $S$ 
of quotients of $\OO_X(-ND)\otimes W$ with Hilbert polynomial $P$.
We have a  universal quotient 
\[
W\otimes p^*\OO_X(-ND)\rightarrow \Qu\rightarrow 0.
\]
As $Q\times X$ is quasi-compact and $G$ reductive, the algebra $\Sym^*(V\otimes \Qu)^G$
is generated by elements of degree at most $k$, for some $k$.
Given $q\in Q$, it then follows that a reduction for the
corresponding vector bundle
\[
\sigma:\Sym^*(V\otimes\Qu|_{q\times_S X})^G\rightarrow \OO_X
\]
is determined by a section of a finite dimensional affine space 
 $$\Sigma = \sspec( \Sym^*_{\OO_S}(\bigoplus_{i=0}^k
       \shom(\Sym^i(V\otimes W)^G\otimes\OO_S, f_*(\OO(iND)) )^\vee
 )).
 $$
The set of all possible such pairs $(\sigma, q)\in \Sigma \times
Q$ coming from an algebra homomorphism
\[
\sigma:\Sym^*(V\otimes\Qu|_{q\times_S X})^G\rightarrow \OO_X
\]
forms a closed   subscheme $\Lambda$. The action of the 
algebraic group $K$ lifts to the Quot scheme and preserves
the subscheme $\Lambda$ and this action also commutes
with the ${\rm GL(W)}$ action.
The subset $U\subset \Lambda$ parameterizing stable bundles is open and 
corresponds to the stable vector bundle locus in a linearization
of the action of ${\rm GL}(W)$.  For this we again
refer the reader to \cite[pg. 1199]{schmitt:02}.
In particular, we can form both the GIT quotient ${\rmbun}_G = U/\!/G$
and the stack theoretic quotient $\bun_G^s = [U/G]$.
The above discussion shows that the natural morphism 
$$\bun_G^s \to {\rmbun}_G $$
is $K$-equivariant.
\end{proof}

We will call the quotient stack $[\rmbun_G/K]$
the coarse moduli space over $[S/K]$. In particular,
this construction  yields the coarse moduli space $\bun_{G,\sC_g}/\M_g$
over the universal curve.

 Let $d$ be the codimension
of the closed complement of $\bun^{s}_G$ in $\bun_G$.
For the statement below, we can either take rational Betti cohomology
with respect to some embedding $\kappa\subset \C$, or \'{e}tale cohomology
with $\Q_\ell$-coefficients.

\begin{prop}
 Then
 \begin{enumerate}
 \item[(i)] The natural map $\bun^{s}_G\rightarrow \rmbun_G$
   induces an isomorphism on rational cohomology in all degrees.
\item[(ii)]
The inclusion $\bun^{s}_G\rightarrow \bun_G$ induces an 
isomorphism on rational cohomology in degrees smaller than $i<2d$.
\end{enumerate}
\end{prop}

\begin{proof}
The  fibre of  map  $\bun^{s}_G\rightarrow \rmbun_G$ over a stable
bundle $P$, can be identified with $B \Gamma(Aut(P))$. Where the group
$\Gamma(Aut(P))$ of global automorphisms of $P$ is finite \cite[prop
3.2]{ramanathan1}. Since the fibres have no higher rational
cohomology,  the first assertion follows by the Leray spectral sequence.
The  second assertion can be deduced from the Gysin sequence.
\end{proof}

Note that $d$ has been calculated in proposition \ref{prop:dforABCD}
for the simple families $A_n,B_n,C_n$ and $D_n$.
We will continue using the notations introduced in subsection 4.6.

\begin{cor}
Let $i<2d$. Then we have an isomorphism of motives
$$h^i(\bun_{G,\sC_g}/\M_g)\cong h^i(\rmbun_{G,\sC_g}/\M_g)$$
\end{cor}

\begin{cor}
Suppose $i<2d$, then
\begin{enumerate}
\item  The motive $h^i(\rmbun_{G,\sC_g}/\M_g)$ is contained in the 
Tannakian subcategory generated by $\pi:\sC_g\to \M_g$.
\item When $\kappa=\C$, the variation of mixed Hodge structure associated to the $i$th
  cohomology of $\rmbun_{G,\sC_g}/\M_g$ is pure and lies in the
  Tannakian subcategory generated by $R^1\pi_*\Q$
\end{enumerate}
\end{cor}

A similar conclusion can be made about the motives of the moduli space
of vector bundles of {\em coprime} rank and degree for all
$i$  \cite{arapura}.  However, the above corollary cannot be extended to all $i$, 
because Cappell, Lee and Miller \cite{cappell} have shown that the
Torelli group would act nontrivially on $H^i(\rmbun_{SL_2})$ for  most $i$
not less than $\dim \rmbun_{SL_2}$. In this situation, we expect the 
variation to be genuinely mixed. In fact, we make the following
conjecture:

\begin{conj}
  The variation of Hodge structure associated to
  $Gr_WH^i(\rmbun_{G,C})$ lies in  the Tannakian subcategory generated by
  $R^1\pi_*\Q$ for all $i$. In particular, the Torelli group acts trivially on this space.
\end{conj}

\end{document}